\documentclass[12pt]{article}

\usepackage[centertags,sumlimits]{amsmath}
\usepackage{amssymb}
\usepackage{amsthm}

\newtheorem{thm}{Theorem}[section]       \newtheorem{lem}[thm]{Lemma}
\newtheorem{prop}[thm]{Proposition}      \newtheorem{cor}[thm]{Corollary}
           \newtheorem{remark}[thm]{Remark}
\newtheorem*{remark*}{Remark}

    \newcommand{\cb}{{\cal B}}   
\newcommand{\ce}{{\cal E}}       
       
       \newcommand{\cp}{{\cal P}}
\newcommand{\calr}{{\cal R}}     \newcommand{\ct}{{\cal T}}
\newcommand{\cu}{{\cal U}}       \newcommand{\cv}{{\cal V}}
\newcommand{\cw}{{\cal W}}    

\newcommand{\R}{\mathbb{R}}     
\newcommand{\al}{{\alpha}}    

  \newcommand{\exv}{\mbox{\sf Exc}_{\cal V}}
\newcommand{\exub}{ \mbox{\sf Exc}_{{\cal U}_{\beta}} }

\newcommand{\ex}{\mbox{\sf Exc}_{\cal U}}
\newcommand{\dis}{\mbox{\sf Diss}_{\cal U}}
\newcommand{\con}{\mbox{\sf Con}_{\cal U}}

\numberwithin{equation}{section}

\begin{document}

\noindent
{\Large\bf  Quasi-regular topologies for $L^p$-resolvents\\[2mm]
and semi-Dirichlet forms}\\

\noindent
LUCIAN BEZNEA$^1$, NICU BOBOC$^2$ and MICHAEL R\"OCKNER$^3$\\
$^1${\it\small 
Institute of Mathematics "Simion Stoilow"
of the Romanian Academy, P.O. Box \mbox{1-764,}  RO-014700 Bucharest, Romania.
(e-mail: lucian.beznea@imar.ro)}\\
$^2${\it\small 
Faculty  of  Mathematics and Informatics, University of Bucharest,
str. Academiei 14, RO-010014  Bucharest, Romania}\\
$^3${\it\small
Fakult\"at f\"ur Mathematik, Universit\"at Bielefeld, 
Postfach 100 131, D-33501 Bielefeld, Germany.
(e-mail: roeckner@mathematik.uni-bielefeld.de)}

\vspace{7mm}

\noindent{\small {\bf Abstract.} We prove that for any
semi-Dirichlet form $( \mbox{\Large $\varepsilon$}, D(\mbox{\Large
$\varepsilon$}) )$ on a measurable Lusin space  $E$ there exists a
Lusin topology  with the given $\sigma$-algebra as the Borel
$\sigma$-algebra so that $( \mbox{\Large $\varepsilon$},
D(\mbox{\Large $\varepsilon$}) )$ becomes quasi--regular. However
one has to enlarge $E$ by a zero set. More generally a
corresponding result for arbitrary $L^p$-resolvents is proven.\\

}\\


\section*{Introduction}
Let $E$ be a Lusin topological space (i.e. the continuous
one-to-one image of a Polish space) with Borel $\sigma$-algebra
${\cal B}$. Let $m$ be a $\sigma$-finite measure on $(E,{\cal B})$
and $L^p(E,m)$, $p\in[1,\infty]$, the corresponding (real)
$L^p$--spaces. Let $( \mbox{\Large $\varepsilon$}, D(\mbox{\Large
$\varepsilon$}) )$  be a semi-Dirichlet form on $L^2(E,m)$ in the
sense of \cite{MaRo 92}. 
Modifying the main result of \cite{AlMaRo 93}, \cite{MaRo 92}, 
in \cite{MaOvRo 95} 
an analytic characterization of all
semi-Dirichlet forms on $L^2(E,m)$ which are associated with a
nice Markov process (more precisely a so-called $m$-special
standard process) was proved. Such semi-Dirichlet forms are called
quasi-regular. An elaborate theory for such Dirichlet forms has
been developed both for its analytic and probabilistic components
with numerous applications (cf. \cite{MaRo 92}). In particular,
invariance properties under change of topology, more precisely,
the invariance under quasi-isomorphism of the theory was
discovered (cf. \cite{AlMaRo 92}, \cite{ChMaRo 94} and the
Appendix in \cite{FuOsTa 94}) and exploited subsequently (see
e.g., Chap. VI in \cite{MaRo 92}).

A fundamental question, however, remained open, namely whether 
it is enough to have a measurable structure only, in the following sense:
Let $(E,{\cal B})$ be merely a Lusin measurable space (i.e. the
image of a Polish space under a measurable homeomorphism) and $(
\mbox{\Large $\varepsilon$}, D(\mbox{\Large $\varepsilon$}) )$
 a semi-Dirichlet form on $L^2(E,m)$ with $m$ a $\sigma$-finite
measure. Can we find a topology on $E$ with Borel $\sigma$-algebra
equal to the given ${\cal B}$ and making $E$ a Lusin topological
space such that $( \mbox{\Large $\varepsilon$}, D(\mbox{\Large
$\varepsilon$}) )$ is quasi-regular with respect to this topology?
As a consequence one could apply all results on quasi-regular
Dirichlet forms only depending on the measurable structure (such
as measure representations for potentials, spectral analysis,
Beurling-Deny type representations etc.) for Dirichlet forms on
arbitrary Lusin-measurable state spaces.

This question has been addressed in \cite{Do 98}  where a
necessary and sufficient condition on $( \mbox{\Large
$\varepsilon$}, D(\mbox{\Large $\varepsilon$}) )$ and ${\cal B}$
was formulated so that the answer to the above question is
positive. This condition is, however, quite close  to what is
needed in the proof and, therefore, not very useful in
applications (see the example in \cite{Do 98}). The main purpose
of this paper is to show that it is always possible to find a
Lusin topology on $E$ making $( \mbox{\Large $\varepsilon$},
D(\mbox{\Large $\varepsilon$}) )$  quasi-regular, however, one has
to enlarge $E$ by a set of $m$-zero measure (cf. Corollary 3.4
below). Our strategy of proof reveals that such an enlargement is
probably necessary in general, though we cannot formally prove
that.

The organization of this paper is as follows. In Section 2 we
first formulate and prove a corresponding result more generally
for $L^p$-resolvents (cf.\ Theorem 2.2) and apply it subsequently
to semi-Dirichlet forms in Section 3 (see Theorem 3.3 and
Corollary 3.4). Our proof relies heavily on results in \cite{BeBo 04a}, 
in particular the characterization of resolvents of kernels which are
associated to right processes. Therefore, in Section 1 we recall
the most essential notions, and list all relevant results. In
particular, we prove that the above characterization of resolvent
kernels can be generalized to the non--transient case (see Theorem
1.3).

\section{Preliminaries on sub-Markovian resolvents of kernels}
Below we follow the terminology of \cite{BeBo 04a}. Let
$\mathcal{U}=(U_\alpha)_{\alpha>0}$ be a  sub-Markovian resolvent
of kernels on a Lusin measurable space $(E, \mathcal{B})$. Recall
that the resolvent $\cu$ is called {\it proper} provided there
exists a strictly positive function $f\in bp\cb$ such that $Uf\leq
1$, where $U=\sup_{\al >0} U_\al$ is the {\it initial kernel} of
$\cu$; $p\cb$ (resp. $bp\cb$ denotes the set of all positive
numerical (resp. bounded positive) $\cb$-measurable functions on
$E$. If $\beta >0$ then the family $\cu_\beta=(U_{\beta
+\al})_{\al>0}$ is also a sub-Markovian resolvent of kernels on
$(E,\cb )$, having $U_\beta$ as (bounded) initial kernel. Recall
also that a function $s\in p\cb$ is termed $\cu$-{\it supermedian}
if $\al U_\al s\leq s$ for all $\al >0$. A $\cu$-supermedian
function $s$ is named $\cu$-{\it excessive} if in addition
$\sup_{\al >0} \al U_\al s=s$. We denote by $\ce(\cu)$ the set of
all $\cb$-measurable $\cu$-excessive functions on $E$. If $s$ is
$\cu$-supermedian then the function $\widehat{s}$ defined by
$\widehat{s}(x)=\sup_{\al >0} \al U_\al s(x)$, $x\in E$, is
$\cu$-excessive and the set $M=\{ x\in E|$ $s(x)\not=
\widehat{s}(x) \}$ is $\cu$-{\it negligible}, i.e. $U_\al(1_M)=0$
for one (and therefore for all) $\al >0$. We denote by $D_\cu$ the
set of all {\it non-branch} points with respect to $\cu$,
$$
D_\cu= \{ x\in E |\, \inf(s,t)(x)= \widehat{\inf(s,t)}(x)
\mbox{ for all } s,t \in \ce(\cu), \widehat{1}(x)=1 \}.
$$
If $\cu$ is proper then, since $\cb$ is countably  generated, we
have $D_\cu\in \cb$ and the set $E\setminus D_\cu$ is
$\cu$-negligible. Notice that in this case $D_\cu =D_{\cu_\beta}$
for all $\beta >0$.

Let $(E', \cb')$ be a second Lusin measurable space such that
$E\subset E'$, $E\in \cb'$ and $\cb=\cb'|_E$.
For all $\al >0$ define the kernel $U'_\al$ on $(E', \cb')$ by
$$
U'_\al f= 1_E U_\al(f|_E)+\frac{1}{1+\al} 1_{E'\setminus E} f\quad f\in p\cb'.
$$
Then the family $\cu'=(U'_\al)_{\al >0}$
is a sub-Markovian resolvent of kernels on $(E', \cb')$,
called the {\it trivial extension of $\cu$ to} $E'$.
If $\beta >0$ then a function $s\in p\cb'$ will be
$\cu'_\beta$-excessive if and only if $s|_E$ is $\cu_\beta$-excessive.
Particularly we have $D_{\cu'_\beta}=D_{\cu_\beta}\cup (E'\setminus E)$ and:
$\sigma(\ce(\cu_\beta))=\cb$ if and only if $\sigma(\ce(\cu'_\beta))=\cb'$.
If $\cu$ is proper then $\cu'$ is also proper.

Let $M\in \cb$ be such that $U_\al (1_{E\setminus M})=0$ on $M$ for one (and therefore for all) $\al>0$.
Then the family  of kernels $\cu|_M=(U_\al|_M)_{\al >0}$ on $(M,\cb|_M)$ is a sub-Markovian
resolvent of kernels, called the {\it restriction of $\cu$ to} $M$;
the kernel $U_\al|_M$ is defined by
$U_\al|_M(g)=U_\al(\overline{g})|_M$ where $\overline{g}\in p\cb$ and $\overline{g}|_M=g$.

Recall that a $\sigma$-finite measure $\xi$ on $(E,\cb)$ is called
{\it $\cu$-excessive} if $\xi \circ \al U_\al \leq \xi$ for all
$\al>0$. We denote by $\ex$ the set of all $\cu$-excessive
measures. Let further $L: \ex\times \ce(\cu)\longrightarrow  \R$
be the {\it energy functional} (associated with $\cu$), $L(\xi,
s)=\sup \{ \mu(s)|$ $\mu$ a $\sigma$-finite measure, $\mu\circ
U\leq \xi \},$ for all $\xi\in \ex$ and $s\in \ce(\cu)$. A
$\cu$-excessive measure  of the form $\mu\circ U$ (where $\mu$ is
a $\sigma$-finite  measure) is called {\it potential}.

For the rest of this section we suppose  that $D_{\cu_\beta}=E$
and $\sigma(\ce(\cu_\beta))=\cb$ for one (and therefore for all)
$\beta>0$.
\\

\noindent
{\bf The transient case}

Suppose that $\cu$ is proper.
Notice  that if $\mu\circ U=\nu\circ U\in \ex$ then $\mu=\nu$.
Moreover the set $\ex$ is an $H$-cone with respect to the usual order relation
on the positive $\sigma$-finite measures; see e.g. \cite{Ge 90}.

A $\cu$-excessive measure $\xi$ is called {\it purely excessive}
(resp. {\it invariant}) if $\inf_\al \xi\circ \al U_\al =0$ (resp.
$\xi\circ \al U_\al =\xi$ for all $\al>0$). Note that if $\xi\in
\ex$ then the measure $\xi_o=\inf_\al \xi\circ \al U_\al$ is
invariant and $\xi -\xi_o$ is purely excessive. Also, every
potential is purely excessive.

The proof of the following lemma is given in the Appendix.

\begin{lem} 
If $\beta>0$ then the following assertions hold.

$a)$ Let $\xi\in \ex$. Then the measure $\xi'=\xi - \xi\circ \beta
U_\beta$ is $\cu_\beta$-excessive. If in addition $\xi$ is purely
excessive then $\xi= \xi'\circ (I+\beta U)$ and for every
$\eta\in\ex$ with $\xi - \xi\circ \beta U_\beta\leq$ $\eta -
\eta\circ \beta U_\beta$ we have $\xi\leq \eta$.

$b)$  If $\xi'\in \exub$ and the measure $\xi= \xi'\circ (I+\beta
U)$ is $\sigma$-finite, then $\xi\in \ex$. Furthermore  it is
purely excessive and $\xi'=\xi - \xi\circ \beta U_\beta$.
\end{lem}

We collect now some results on the semisaturation and saturation of $E$; cf. \cite{BeBo 04a}.
The set $E$ is called {\it semisaturated} (resp. {\it saturated}) {\it with respect to}
$\cu$ provided that every $\cu$-excessive measure dominated by a potential is also a potential
(resp. every $\xi\in \ex$ with $L(\xi, 1) < \infty$ is a potential).
If $\xi\in \ex$ then $E$ is termed {\it $\xi$-semisaturated} if every $\cu$-excessive measure
dominated by a potential dominated by $\xi$ is also a potential.
The following assertions hold.

$i)$ If $E$ is saturated with respect to $\cu$ then $E$ is semisaturated with respect to $\cu$.

$ii)$ The set $E$ is semisaturated with respect $\cu$ if and only
if there exists a Lusin topology on $E$ such that $\cb$ is the
$\sigma$-algebra of all Borel sets on $E$ and there exists a right
process with state space $E$, having $\cu$ as associated
resolvent.

$iii)$  There exist a second Lusin measurable  space $(E_1,\cb_1)$
such that $E\subset E_1$, $E\in \cb_1$, $\cb=\cb_1|_{E}$,  and a
proper sub-Markovian resolvent of kernels
$\cu^1=(U^1_\al)_{\al>0}$ on $(E_1,\cb_1)$ such that
$D_{\cu^1}=E_1$, $\sigma(\ce(\cu^1))=\cb_1$,
$U^1_\al(1_{E_1\setminus E})=0$, $E_1$ is saturated with respect
to $\cu^1$ and $\cu$ is the restriction of $\cu^1$ to $E$. 
In particular, $\cu^1$ is the resolvent of a right process with state
space $E_1$ for a suitable Lusin topology on $E_1$. More precisely
one can take $E_1$ as the set of all extreme points of the set $\{
\xi\in \ex| L(\xi ,1) \leq 1 \}$, endowed with the
$\sigma$-algebra $\cb_1$ generated by the functionals
$\widetilde{s}$, $\widetilde{s}(\xi)=L(\xi, s)$ for all $\xi\in
E_1$ and $s\in \ce(\cu)$. The set $E_1$ is called the {\it
saturation} of $E$.

$iv)$ Let $(E',\cb')$ be a Lusin measurable space such that
$E\subset E'$, $E\in \cb'$, $\cb=\cb'|_{E}$,  and there exists a
proper sub-Markovian resolvent of kernels $\cu'=(U'_\al)_{\al>0}$
on $(E',\cb')$ with  $D_{\cu'}=E'$, $\sigma(\ce(\cu'))=\cb'$,
$U'_\al(1_{E'\setminus E})=0$, $E'$ is saturated with respect to
$\cu'$ and $\cu$ is the restriction of $\cu'$ to $E$. Then the map
$x\longmapsto \varepsilon_x \circ U'$ is a measurable  isomorphism
between $(E', \cb')$ and the measurable space  $(E_1, \cb_1)$
defined in $iii)$ above.

$v)$ The set $E$ is semisaturated (resp. $\xi$-semisaturated,
where $\xi$ is a fixed $\cu$-excessive measure) if and only if $E_1\setminus E$
is a polar (resp. $\xi$-polar) subset of $E_1$ (with respect to $\cu^1$);
recall that a set $M\in \cb$ is {\it polar} (resp. {\it $\xi$-polar}) with respect to $\cu$
if $\widehat{R^M 1}=0$ (resp. $\widehat{R^M 1}=0$ $\xi$-a.e.), where $R^M 1$ denotes the reduced function
(with respect to $\cu$) of $1$ on $M$,
$R^M1=\inf \{ s\in \ce(\cu)|$ $s\geq 1$ on $M \}.$

$vi)$ If $E$ is $\xi$-semisaturated then there exists a proper sub-Markovian resolvent of kernels
$\cu'=(U'_\al)_{\al>0}$ on $(E,\cb)$ such that the set $E$ is semisaturated with respect to $\cu'$
and for all $f\in p\cb$ and $\al>0$ the set $[U_\al f\not= U'_\al f]$ is $\xi$-polar.

$vii)$ Let $A\in \cb$ be such that $U_\al(1_{E\setminus A})=0$ on $A$ and $\cu'$
the trivial extension of $\cu|_A$ to $E$.
Then $A$ is semisaturated with respect to $\cu|_A$ if and only if
$E$ is semisaturated with respect to $\cu'$.

\begin{prop} 
Let $\beta>0$. Then $E$ is semisaturated (resp. saturated) with
respect to $\cu$ if and only if it is semisaturated (resp.
saturated) with respect to $\cu_\beta$.
\end{prop}

\begin{proof}
Suppose that $E$ is semisaturated with respect to $\cu$ and let
$\xi', \mu\circ U_\beta \in \exub$, $\xi'\leq \mu\circ U_\beta$.
Clearly we may assume that $\mu$ is finite and thus $\xi'$ is also a finite measure.
By Lemma 1.1 it follows that that the measure
$\xi=\xi'\circ(I+\beta U)$ is $\cu$-excessive and $\xi'=\xi- \xi\circ \beta U_\beta$.
Since $\xi\leq \mu\circ U_\beta (I+\beta U)=$ $\mu \circ U$ we deduce by hypothesis
that there exists a $\sigma$-finite measure $\nu$ on $(E,\cb)$  such that $\xi=\nu\circ U$
and thus $\xi'=\nu\circ U(I-\beta U_\beta)=$ $\nu\circ U_\beta.$

If $E$ is saturated with respect to $\cu$ and $\xi'\in \exub$ is
such that $L_\beta(\xi', 1)<\infty$ ($L_\beta$ denotes the energy
functional associated with $\cu_\beta$) then we claim that the
measure $\xi=\xi'\circ (I+\beta U)$ is $\sigma$-finite. Indeed,
let $(\mu_n)_n$ be a sequence of positive measures on $(E,\cb)$
such that $\mu_n \circ U_\beta \nearrow \xi'.$ From
$\mu_n(1)=L_\beta(\mu_n\circ U_\beta , 1)\leq L_\beta(\xi', 1)$ it
follows that $\sup_n \mu_n(1)<\infty$. If $f_o\in bp\cb$ is such
that $Uf_o\leq 1$ then we get $\xi(f_o)=\xi'\circ(I+\beta
U)(f_o)=$ $\sup_n \mu_n\circ U_\beta (I+\beta U)(f_o)=$ $\sup_n
\mu_n(Uf_o)\leq$ $\sup_n \mu_n(1)< \infty$. Hence the measure
$\xi$ is $\sigma$-finite and by Lemma 1.1 we obtain that $\xi$ is
$\cu$-excessive and $\xi'= \xi \circ (I- \beta U_\beta)$. Since
$L(\xi, 1)=\sup_n \mu_n(1)<\infty$ and   $E$ is saturated with
respect to $\cu$, it follows that there exists a $\sigma$-finite
measure $\mu$ on $(E,\cb)$ such that $\xi=\mu\circ U$ and thus
$\xi'=\mu\circ U_\beta$.

Assume now that $E$ is semisaturated with respect to $\cu_\beta$
and let $\xi, \mu\circ U \in \ex$, $\xi\leq \mu\circ U$. The
measure $\xi$ is purely excessive and we may suppose that $\mu$ is
finite. Consequently the measure $\mu'=\mu\circ (I+\beta U)$ is
$\sigma$-finite. Again by Lemma 1.1 it follows that the measure
$\xi'= \xi \circ (I- \beta U_\beta)$ is $\cu_\beta$-excessive.
Since $\xi'\leq \mu\circ U=\mu'\circ U_\beta$, by hypothesis there
exists a $\sigma$-finite measure $\nu$ on $(E,\cb)$  such that
$\xi'=\nu\circ U_\beta.$ As a consequence we get $\xi=\xi'\circ
(I+\beta U)=\nu\circ U.$

Let us suppose now that $E$ is saturated with respect to
$\cu_\beta$ and $\xi\in \ex$ is such that $L(\xi, 1)<\infty$. If
$E_1$ is the saturation of $E$ with respect to $\cu$ then $\xi$ is
a potential on $E_1$ and thus it is purely excessive. Lemma 1.1
implies that the measure $\xi'=\xi-\xi\circ \beta U_\beta$ belongs
to $\exub$ and $\xi=\xi' \circ (I+\beta U)$. We consider a
sequence $(\mu_n)_n$ of positive $\sigma$-finite measures on
$(E,\cb)$ such that $\mu_n \circ U_\beta \nearrow \xi'.$
Consequently, we have $\mu_n \circ U \nearrow \xi$ and
$L_\beta(\xi',1)=\sup_n \mu_n(1)=L(\xi,1)<\infty$. Therefore,
there exists a $\sigma$-finite measure $\mu$ on $(E,\cb)$ such
that $\xi'=\mu\circ U_\beta$ and so $\xi= \xi' \circ (I+\beta
U)=\mu\circ U$.
\end{proof}

\noindent {\bf The non-transient case}

Firstly recall  some facts
on  Ray cones. Assume that the initial kernel $U$ of $\cu$ is
bounded. A {\it Ray cone associated with $\cu$} is a convex cone
$\calr$ of bounded $\cu$-excessive functions such that: $U_\al
(\calr)\subset \calr$ for all $\al>0$, $U((\calr -\calr)_+)\subset
\calr,$ $\sigma(\calr)=\cb$, $\calr$ is min-stable, separable in
the uniform norm and contains the positive constant functions.

We state here a slightly modified version of Proposition 1.5.1 in \cite{BeBo 04a}:
Let $\beta >0$. Then there  exists a Ray cone $\calr_\beta$ associated with $\cu_\beta$,
such that $U_\al(\calr_\beta)\subset \calr_\beta$ for all $\al>0$.

We claim that the above assertion $ii)$ is true without assuming
that $\cu$ is proper. Namely the following result is a variant of
assertion $ii)$, in the case when the initial kernel $U$ is not
necessary a proper one; compare with \cite{St 89}.

\begin{thm} 
The set  $E$ is semisaturated with respect to $\cu_\beta$ if and
only if there exists a Lusin topology on $E$ such that $\cb$ is
the $\sigma$-algebra of all Borel sets on $E$ and there exists a
right process with state space $E$, having $\cu$ as associated
resolvent.
\end{thm}

\begin{proof}
It is known that $E$ is semisaturated with respect to $\cu_\beta$ whenever
$\cu$ is the resolvent of a right process; see \cite{Ge 90}.
For the converse statement we shall adapt the proofs of
Theorem 1.8.11 and Corollary 1.8.12 in \cite{BeBo 04a}.

First assume  that $E$ is saturated with respect to $\cu_\beta$.
Let $\calr_\beta$ be a Ray cone associated with $\cu_\beta$  such
that $U_\al(\calr_\beta)\subset \calr_\beta$ for all $\al>0$, and
$Y$ the (Ray) compactification of $E$ with  respect to
$\calr_\beta$. By Proposition 1.5.8 in \cite{BeBo 04a} there
exists a Ray resolvent $\widetilde{\cu}= (\widetilde{U}_\al)_{\al
>0}$ on $Y$ such that
$\widetilde{U}_\al(\widetilde{s})=\widetilde{U_\al s}$ for all
$s\in \calr_\beta$ and $\al>0$, where for each $s\in \calr_\beta$
we have denoted by $\widetilde{s}$ the unique continuous extension
of $s$ to $Y$. Particularly $\widetilde{U}_\al(1_{Y\setminus
E})=0$ on $E$ for all $\al>0$ and $\cu$ is the restriction to $E$
of $\widetilde{U}$. Consequently (see e.g. \cite{Sh 88}) the
restriction of $\widetilde{\cu}$ to $D=D_{\widetilde{\cu}}$ is the
resolvent of a right process $X$ with state space $D$, endowed
with the Ray topology induced by $\calr_\beta$ (i.e. the trace on
$D$ of the topology on $Y$). From Theorem 1.8.11 in \cite{BeBo
04a} we have $E=\{ x\in D|$ $\widetilde{U}_\al(1_{D\setminus
E})(x) =0 \}$. In addition $E$ is a Borel subset of $Y$,
$\widetilde{U}_\al(1_{D\setminus E})=0$ on $E$ and it is a finely
closed set with respect to $\widetilde{\cu}_\beta$; the {\it fine
topology} is the topology generated by
$\ce(\widetilde{\cu}_\beta)$. As a consequence we may consider the
restriction of $X$ to $E$ and $\cu$ becomes the resolvent of this
right process, since $\widetilde{\cu}|_E=\cu.$

If $E$ is only semisaturated with respect to $\cu_\beta$, then we consider the saturation $E_1$ of $E$
with respect to $\cu_\beta$ and let
$\cu^1=(U^1_\al)_{\al>0}$ be the resolvent of kernels on on $(E_1,\cb_1)$ such that $D_{\cu^1}=E_1$,
$\sigma(\ce(\cu^1_\beta))=\cb_1$, $U^1_\al(1_{E_1\setminus E})=0$ and $\cu^1|_E=\cu$.
By the first part of the proof there exists a right process $X$ with state space $E_1$
(endowed with a Ray topology), having $\cu^1$ as associated resolvent.
By $v)$ we deduce that the set $E_1\setminus E$ is polar (with respect to $\cu^1_\beta$)
and therefore the restriction of $X$ to $E$ is a right process with state space $E$ and having $\cu$
as associated resolvent, completing the proof.
\end{proof}

\begin{remark*}{\rm
By Proposition 1.2 it follows that the condition of semisaturation with respect to $\cu_\beta$
in Theorem 1.3 does not depend on $\beta$.
}
\end{remark*}

Recall that a $\cu$-excessive measure $\xi$ is called {\it
dissipative} (resp. {\it conservative}) provided that $\xi=\sup \{
\mu\circ U| $ $\ex \ni\mu \circ U \leq \xi \}$ (resp. there is no
non-zero potential  $\cu$-excessive measure dominated by $\xi$).
The set of all dissipative (resp. conservative) $\cu$-excessive
measures is denoted by $\dis$ (resp. $\con$). As in \cite{Ge 90}
one can show that $\dis$ and $\con$ are solid convex subcones of
$\ex$, $\dis\cap \con =0$ and every $\xi\in \ex$ has a unique
decomposition of the form $\xi=\xi_d+\xi_c$, where $\xi_d\in\dis$ and $\xi_c\in\con$. 
Moreover, if $f\in p\cb$ is strictly positive
and  $\xi(f)<\infty$ then $\xi_d=\xi|_{[Uf<\infty]}$ and
$\xi_c=\xi|_{[Uf= \infty]}$; See also Proposition A1 in the
Appendix.

The next result is an extension of assertion $vi)$ to the
non-transient case.

\begin{prop} 
Let $\xi\in \dis$ be such that $E$ is $\xi$-semisaturated with respect to $\cu$
(i.e., every $\cu$-excessive measure dominated by a potential
dominated by $\xi$ is also a potential).
Then there exists a proper sub-Markovian resolvent of kernels $\cu'=(U'_\al)_{\al >0}$ on $(E,\cb)$
such that $E$ is semisaturated with respect to  $\cu'$ and the set
$[U_\al f\not= U'_\al f]$
is $\xi$-polar with respect to $\cu_\beta$ for all $f\in p\cb$ and $\al >0$.
Moreover there exists a $\xi$-polar finely closed set $A\in \cb$
such that $U(1_A)=0$ on $E\setminus A$ and
$\cu'$ may be chosen as the trivial extension to $E$ of the restriction of $\cu$ to $E\setminus A$.
\end{prop}

\begin{proof}
Let $f\in p\cb$ be strictly positive such that $\xi(f) < \infty$.
The set $A=[Uf= \infty]$ is finely closed, $U(1_A)=0$ on
$E\setminus A$ and from $\xi\in \dis$ we get  $\xi(A)=0$.
Therefore the set  $A$ is $\xi$-polar with respect to $\cu_\beta$.
If  $\cv$ is  the restriction of $\cu$ to $E\setminus A$ then we
deduce that $\cv$ is a proper sub-Markovian resolvent of kernels
on $(E\setminus A, \cb|_{E\setminus A})$ such that
$\sigma(\ce(\cv))= \cb|_{E\setminus A}$ and $D_\cv= E\setminus A$.
Clearly the measure $\xi$ belongs to $\exv$. We show that
$E\setminus A$ is $\xi$-semisaturated with respect to $\cv$.
Indeed, let $\eta, \mu\circ V\in \exv$, with $\eta\leq \mu\circ
V\leq \xi$, where $\mu$ is a $\sigma$-finite measure on
$E\setminus A$. We deduce that  $\eta, \mu\circ U\in \ex$ and
$\eta\leq \mu\circ U\leq \xi$. Since $E$ is $\xi$-semisaturated
with respect to $\cu$, there exists a $\sigma$-finite measure
$\nu$ on $E$ such that $\eta=\nu\circ U$. Since the  set $A$ is
$\mu$-polar and $\mu$-negligible, it follows that it is also
$\nu$-negligible and  consequently $\eta=\nu|_{E\setminus A} \circ V$. 
y $vi)$ there exists a proper sub-Markovian resolvent of
kernels $\cv'=(V'_\al)_{\al>0}$ on $( E\setminus
A,\cb|_{E\setminus A} )$ such that $E\setminus A$ is semisaturated
with respect to $\cv'$ and $[V_\al f = V'_\al f]$ on $E_o$ for all
$f\in p\cb|_{E\setminus A}$ and $\al>0$, where $E_o\in \cb$ is
such that $E_o \subset E\setminus A$, $E\setminus E_o$  is
$\xi$-polar and $U(1_{E\setminus E_o})=0$ on $E_o$. From $vii)$ we
conclude  that  the trivial extension $\cu'$ of $\cv'|_{E_o}$ to
$E$ satisfies the required conditions.
\end{proof}

\section{Right processes associated with $L^p$-resolvents}

In the sequel $\mu$ will be a $\sigma$-finite measure on $(E,\cb)$.

Let $\cu'=(U'_\al)_{\al >0}$ be a second sub-Markovian resolvent of kernels on $(E,\cb)$.
We say that $\cu$ and $\cu'$ are {\it $\mu$-equivalent} provided that $U_\al f=U'_\al f$
$\mu$-a.e. for all $f\in p\cb$ and $\al >0$.

\begin{remark*}{\rm
There are examples of two sub-Markovian resolvents of kernels  on the same
space $E$, which are $\xi$-equivalent (where $\xi$ is a $\sigma$-finite measure)
and such that $E$ is semisaturated with respect to only one of them.
Indeed, let $\cu^o$ be a sub-Markovian resolvent on a Lusin measurable space
$(F, \cb_o)$ such that $F$ is not semisaturated with respect to $\cu^o$.
We denote by $E$ the saturation of $F$ with respect to $\cu^o$ (i.e. $E=F_1$)
and let $\cu$ be the resolvent on $E$ such that $\cu|_F=\cu^o$ and $E\setminus F$ is $\cu$-negligible.
Let further $\cu'$ be the trivial extension of $\cu^o$ to $E$.
Then by $vii)$ the set $E$ is not semisaturated with respect to $\cu'$.
Clearly, since $U_\al(1_{E\setminus F})=0$, we deduce that $\cu$ and $\cu'$ are $\xi$-equivalent
with respect to every $\xi\in \ex$.
}
\end{remark*}

\begin{lem} 
Let $N$ be a bounded kernel on  $(E,\cb)$ such that if $B\in \cb$
and  $\mu(B)=0$ then $N(1_B)=0$ $\mu$-a.e.
If $E_o\subset E$, $E_o\in \cb$, is such that $\mu(E\setminus E_o)=0$ then there exists $F\in \cb$,
$F\subset E_o$, such that $\mu(E\setminus F)=0$ and $N(1_{E\setminus F})=0$ on $F$.
\end{lem}

\begin{proof}
Since $\mu(E\setminus E_o)=0$ we get by hypothesis that
$N(1_{E\setminus E_o})=0$ $\mu$-a.e.
Let $(E_n)_{n\geq 1}\subset \cb$ be the sequence defined inductively by
$E_{n+1}=E_n\cap [N(1_{E\setminus E_n})=0]$ if $n\geq 0$.
We have $\mu(E\setminus E_n)=0$  for all $n$ and let $F=\bigcap_{n}E_n$.
Then $F\subset E_o$, $F\in \cb$, $\mu(E\setminus F)=0$ and if $x\in F$ then
$N(1_{E\setminus E_n})(x)=0$ for all $n$.
Therefore $N(1_{E\setminus F})(x)=$ $N(1_{\bigcup_n E\setminus E_n})(x)=$
$\sup_n N(1_{E\setminus E_n})(x)=0$.
\end{proof}

\begin{remark*}{\rm
A procedure similar to Lemma 2.1 has been considered in \cite{Ku 86} and \cite{Mo 93}.
}
\end{remark*}

\begin{thm} 
Let $p\in [1,+\infty]$ and $(V_\al)_{\al>0}$ be a sub-Markovian
strongly continuous resolvent of contractions on $L^p(E,\mu)$,
where $(E,\cb)$ is a Lusin measurable space and $\mu$ is a
$\sigma$-finite measure on $(E,\cb)$. Then there exist a Lusin
topological space $E_1$ with $E\subset E_1$, $E\in \cb_1$ (the
$\sigma$-algebra of all Borel subsets of $E_1$), $\cb=\cb_1|_{E}$,
and a right process with state space  $E_1$ such that its
resolvent, regarded on $L^p(E_1,\overline{\mu}),$ coincides with
$(V_\al)_{\al>0}$, where $\overline{\mu}$ is the measure on $(E_1,
\cb_1)$ extending $\mu$ by zero on $E_1\setminus E$.
\end{thm}

\begin{proof}
Let $(f_k)_k\subset bp\cb\cap L^p(E,\mu)$ be a sequence separating
the points of $E$. For every $\al >0$ we consider a kernel
$\overline{V}_\al$ on $(E,\cb)$ such that $\overline{V}_\al$
coincides with $V_\al$ as an operator on $L^p(E, \mu)$. By
Proposition 1.4.13 in \cite{BeBo 04a} there exists a sub-Markovian
resolvent $\cw=(W_\al)_{\al >0}$ on $(E,\cb)$  such that $W_\al
f=\overline{V}_\al f$ $\mu$-a.e. for all $f\in p\cb$. Let us
consider the set
$$
E_o=\{ x\in E|\, \liminf_n nW_n f_k(x)=f_k(x) \mbox{ for all } k \}.
$$
We have $E_o\in \cb$ and since $(W_\al)_{\al >0}$ is a strongly
continuous resolvent of contractions on $L^p(E,\mu)$ it follows
that $\mu(E\setminus E_o)=0$. Let $\cb'$ be the $\sigma$-algebra
on $E$ generated by $W_1(p\cb)$. Then $\cb'$ is countably
generated and $\cb'|_{E_o}$ separates the points of $E_o$.
Therefore $E_o\in \cb'$ and $\cb'|_{E_o}=\cb|_{E_o}$. By Lemma 2.1
there exists $F\in \cb$, $F\subset E_o$,  such that
$\mu(E\setminus F)=0$ and $W_\al (1_{E\setminus F})=0$ on $F$ for
all $\al >0$. Let $\beta>0$ and $F_o$ be the set of all non-branch
points of $F$ with respect to $\cw_\beta|_F$. Then $F_o\in \cb$,
$\cw_\beta|_F$ is a sub-Markovian resolvent of kernels on $(F_o,
\cb|_{F_o})$, $\ce(\cw_\beta|_F)$ is min-stable, contains the
positive constant functions and generates $\cb|_{F_o}$. Let
$\cu=(U_\al)_{\al >0}$ be the trivial extension of $\cw|_{F_o}$ to
$E$. Then $\cu$ is a sub-Markovian resolvent of kernels on
$(E,\cb)$ such that $D_{\cu_\beta} =E$,
$\sigma(\ce(\cu_\beta))=\cb$ and $(U_\al)_{\al
>0}$ coincides with $(V_\al)_{\al >0}$ as a resolvent on
$L^p(E,\mu)$. We consider now the set $E_1$, i.e.  the saturation
of $E$ with respect $\cu_\beta$ (see $iii)$ in Section 1) and the
resolvent of kernels $\cu^1=(U^1_\al)_{\al >0}$ on $(E_1,\cb_1)$
whose restriction to $E$ is $\cu$ and $U^1_\al (1_{E_1\setminus
E})=0$. Since $E_1$ is saturated with respect to $\cu^1_\beta$, we
deduce from $i)$ and Theorem 1.3 that there exists a Lusin
topology on $E_1$ such that $\cb_1$ is the $\sigma$-algebra of all
Borel sets on $E_1$ and $\cu^1$ is the resolvent of a right
process with state space $E_1$. Clearly $U^1_\al =V_\al$ for all
$\al >0$, regarded as an equality of operators on $L^p(E_1,
\overline{\mu}).$
\end{proof}

\begin{remark}{\rm
Under the assumptions of Theorem 2.2 we have proved that there
exists a sub-Markovian resolvent of kernels $\cu=(U_\al)_{\al >0}$
on $(E,\cb )$ such that for $\beta >0$ we have $D_{\cu_\beta}=E$,
$\sigma(\ce(\cu_\beta))=\cb$ and $U_\al=V_\al$ as operators on
$L^p(E,\mu)$ for all $\al>0$. Moreover the following assertions
hold.

$a)$ $\cu=(U_\al)_{\al >0}$ is the resolvent of a right process with state space $E$ if and only if
$E$ is semisaturated with respect to $\cu_\beta$ (cf. Theorem 1.3).

$b)$ If $\mu$ is $\cu_\beta$-excessive and $E$ is $\mu$-semisaturated with respect to $\cu_\beta$
(or if $\mu\in \dis$ and $E$ is $\mu$-semisaturated with respect to $\cu$)
then by $vi)$, Proposition 1.4  and Theorem 1.3
there exist a Lusin topology on $E$ and a right process with state space
$E$ such that its resolvent and $\cu$ are $\mu$-equivalent.
}
\end{remark}

The following result is a consequence of Proposition 7.5.2 in \cite{BeBo 04a},
Theorem 2.2 and Remark 2.3.

\begin{cor} 
Let $\cu=(U_\al)_{\al >0}$ be a sub-Markovian resolvent of kernels
on $(E,\cb )$ such that for $\beta >0$ we have $D_{\cu_\beta}=E$ and $\sigma(\ce(\cu_\beta))=\cb$.
If $\mu\in \ex$  then there exists a second sub-Markovian resolvent of kernels
$\cu^*=(U_\al^*)_{\al >0}$ on $(E,\cb )$ such that for $\beta >0$ we have $D_{\cu^*_\beta}=E$,
$\sigma(\ce(\cu^*_\beta))=\cb$ and  $\int_E fU_\al g d\mu =$ $\int_E gU^*_\al f d\mu$
for all $f,g \in  p\cb$ and $\al >0$.
\end{cor}

\section{Tightness of capacity and quasi-regularity}
In this section we shall give   conditions  on an $L^p$-resolvent
to ensure tightness of the capacity induced by the reduction
operator, the existence  of  quasi-continuous versions for the
elements being in the domain of the generator and the standardness
property of the associated right process.

Let $(V_\al)_{\al>0}$ be a sub-Markovian resolvent on $L^p(E,\mu)$
as in Theorem 2.2 and $\beta>0$. An element $u\in L^p_+(E, \mu)$
is called a $\beta$-{\it potential} provided that $\al
V_{\beta+\al} u\leq u$ for all $\al>0$. We denote by $\cp_\beta$
the set of all, $\beta$-potentials. It is known that (see e.g.
Proposition 3.1.10 in \cite{BeBo 04a}) the ordered convex cone
$\cp_\beta$ is a cone of potentials in the sense of G. Mokobodzki,
cf. \cite{Mo 70} (see also \cite{BeBo 04a}). Particularly if $u,
u'\in \cp_\beta$, $u\leq u'$, then there exists $R_\beta (u-u')\in
\cp_\beta$, i.e. the {\it r\'eduite} of $u-u'$, defined by
$R_\beta (u-u')=\bigwedge \{ v\in \cp_\beta|$ $v\geq u-u'\}$ (here
$\bigwedge$ denotes the infimum in $\cp_\beta$). An element $u\in
\cp_\beta$ is called {\it regular} if for every sequence
$(u_n)_n\subset \cp_\beta$ with $u_n\nearrow u$ we have
$R_\beta(u-u_n)\searrow 0$.

\begin{remark}{\rm  
$a)$ If $f\in L^p(E,\mu )$ then $V_\beta f$ is regular.
If $u\in \cp_\beta$ then $V_\al u$ is regular for every $\al>0$.

$b)$ Let $u\in \cp_\beta$. If there exists a sequence $(u_n)_n$ of
regular elements from $\cp_\beta$ with $u_n\nearrow u$ and
$R_\beta(u-u_n)\searrow 0$ then by Proposition 3.2.3 in \cite{BeBo
04a} it follows that $u$ is regular. Consequently by  $a)$ we
deduce that: $u$ is regular if and only if $R_\beta(u-nV_n
u)\searrow 0$.

$c)$ Assume that $\cv_\beta=(V_{\beta+\al})_{\al >0}$
is the resolvent of a right process and let
$u \in \ce(\cv_\beta)\cap L^p(E,\mu)$, $u < \infty$.
Then $u$ is regular if and only if there exists
a continuous additive functional whose potential equals $u$ $\mu$-a.e.
}
\end{remark}

Let $f_o\in L^p(E,\mu)$ be strictly positive.
We consider the following property of the resolvent $(V_\al)_{\al>0}$:
\begin{center}
$(*)$\hspace{5mm} {\it every $\beta$-potential dominated by $V_\beta f_o$ is regular.}
\end{center}

\begin{remark*}{\rm
Since $V_\beta f_o >0$ it follows from \cite{BeBo 04a} that condition $(*)$ is equivalent with the
following one: every $\beta$-potential dominated by a regular element from $\cp_\beta$ is also regular.
}
\end{remark*}

\begin{prop} 
Condition $(*)$ does not depend on $\beta.$
\end{prop}

\begin{proof}
Let $\beta'>\beta>0$ and assume that condition $(*)$ holds for $\beta$.
If $(u_n)_n \subset \cp_{\beta'}$, $u_n\nearrow u\in \cp_{\beta'}$, $u\leq V_{\beta'}f_o$,
then the element $v=u+(\beta'-\beta)V_\beta u$ belongs to $\cp_\beta$, $v\leq V_\beta f_o$
and thus $v$ is regular.
Setting $v_n=u_n+(\beta'-\beta)V_\beta u$ we get $v_n\nearrow v$,
$v_n= u_n +(\beta'-\beta)V_\beta u_n + (\beta'-\beta)V_\beta(u-u_n)\in \cp_\beta$ and since
$\cp_\beta\subset \cp_{\beta'}$ it follows that
$R_{\beta'}(u-u_n)=R_{\beta'}(v-v_n)\leq$ $R_\beta(v-v_n)\searrow 0$.

Assume now that condition $(*)$ holds for $\beta'$ and let
 $(u_n)_n \subset \cp_{\beta}$, $u_n\nearrow u\in \cp_{\beta}$.
Then the element $v=u-(\beta'-\beta)V_{\beta'}u$ belongs to
$\cp_{\beta'}$. If  $u\leq V_{\beta}f_o$, since by Remark 3.1
$V_\beta f_o$ is a regular element of $\cp_{\beta'}$, we deduce
that $v$ is regular in $\cp_{\beta'}$. Let $(f_n)_n\subset
L^p(E,\mu)$ be such that $V_{\beta'}f_n\nearrow v$. Then
$R_{\beta'}(v-V_{\beta'} f_n)\searrow 0$ and $V_\beta f_n \nearrow
u$. To show that $u$ is regular, again by Remark 3.1 it suffices
to prove that $R_{\beta}(u-V_{\beta} f_n)\searrow 0$. Notice that
if $u', u''\in \cp_\beta$, $f=u'-u''$, then $R_\beta(f)\leq
(I+(\beta'-\beta)V_\beta)R_{\beta'}(f-(\beta'-\beta)V_{\beta'}
f).$ We conclude that $R_\beta(u-V_\beta f_n)\leq
(I+(\beta'-\beta)V_\beta)R_{\beta'}(v-V_{\beta'} f_n)\searrow 0.$
\end{proof}

\begin{remark*}
{\rm
$a)$ Let $\cu$   be the resolvent of kernels from Remark 2.3,
$\cu^*$ be a second resolvent given by Corollary 2.4
and suppose that they are associated with two right processes  with state space $E$.
Then condition $(*)$ is equivalent with the fact that "the axiom of polarity" holds for
$\cu^*_\beta$, i.e. every semipolar set is $\mu$-polar
with respect to $\cu^*_\beta$ (see \cite{BeBo 04a}).

$b)$ If $(V_\al)_{\al>0}$ is the resolvent of a semi-Dirichlet
form on $L^2(E, \mu)$ then it was shown  in \cite{BeBo 04a} that
condition $(*)$ holds and derived that a semi-Dirichlet form
associated with a right process is quasi-regular; compare with
\cite{Fi 01}, \cite{MaOvRo 95} and \cite{MaRo 92}. }
\end{remark*}

Assume further that in addition $f_o\in L^1(E,\mu)$,
$\lambda_o=f_o\cdot \mu$ and $m=\lambda_o\circ V_\beta$.

The next result is a consequence of Section 3.5 and Theorem 3.7.8
in \cite{BeBo 04a} and Theorem 2.2; see also \cite{BeBo 04b}.

\begin{thm} 
Under the assumptions from Theorem 2.2  suppose that condition
$(*)$ holds. If $(E_1, \ct)$ is the Lusin topological space and
$\cu^1$ the resolvent of the right process with state space $E_1$
given by Theorem 2.2, then the following assertions hold.

$a)$ There exists an increasing sequence $(K_n)_n$ of $\ct$-compact subsets of $E_1$ such that
$$
\inf_n R_\beta^{E_1\setminus K_n} p_o=0\quad (m+\lambda_o)\mbox{-a.e.}
$$
where $p_o=U^1_\beta \widetilde{f}_o$ ($\widetilde{f}_o\in p\cb_1$, $\widetilde{f}_o|_E= f_o$)
and $R^M_\beta p_o$ denotes the reduced function
(with respect to $\cu_\beta^1$) of $p_o$ on the set $M$.

$b)$ Every $\cu^1_\beta$-excessive function $s$ is $\ct$-quasi continuous,
that is there exists an increasing sequence $(K_n)_n$ of $\ct$-compact subsets of $E_1$
such that $s|_{K_n}$ is $\ct$-continuous for all $n$ and
$\inf_n R_\beta^{E_1\setminus K_n} p_o=0$ $(m+\lambda_o)$-a.e.
Particularly, every element from $V_\al(L^p(E,\mu))$
(the domain of the generator of the resolvent  $(V_\al)_{\al>0}$)
possesses a $\ct$-quasi continuous $\mu$-version.

$c)$ The right process having $\cu^1$
as associated resolvent is $(m+\lambda_o)$-special standard.
\end{thm}

As a consequence of the previous theorem and the main result in
\cite{MaOvRo 95} and \cite{MaRo 92} we obtain:

\begin{cor} 
Let $( \mbox{\Large $\varepsilon$}, D(\mbox{\Large $\varepsilon$}))$ 
be a semi-Dirichlet form on $L^2(E,\mu)$, where $\mu$ is a
$\sigma$-finite measure on the Lusin measurable space $(E,\cb)$.
Then there exists a (larger) Lusin topological space $E_1$ such
that  $E\subset E_1$, $E$ belongs to $\cb_1$ (the $\sigma$-algebra
of all Borel subsets of $E_1$), $\cb=\cb_1|_E$, and $(
\mbox{\Large $\varepsilon$}, D(\mbox{\Large $\varepsilon$}) )$
regarded as a semi-Dirichlet form on $L^2(E_1 , \overline{\mu})$
is quasi-regular, where $\overline{\mu}$ is the measure on $(E_1,
\cb_1)$ extending $\mu$ by zero on $E_1\setminus E$.
\end{cor}

\vspace{3mm}

\noindent
{\large\bf Appendix}\\

\noindent
{\bf Proof of Lemma 1.1.}

$a)$ If $\al >0$ then we have $\xi'\circ\al U_{\beta+\al}=$
$\xi\circ \al U_{\beta+\al}- \xi\circ \beta\al U_\beta
U_{\beta+\al}=$ $\xi\circ (\al+\beta)U_{\beta+\al} -\xi\circ \beta
U_\beta\leq$ $\xi -\xi\circ \beta U_\beta=\xi'.$ 
For $\alpha<\beta$  we have also
$\xi'\circ(I+(\beta-\al)U_\al)=$ $\xi\circ (I-\beta
U_\beta+(\beta-\al)U_\al- (\beta-\al)\beta U_\al U_\beta)=$ $\xi-
\xi\circ\al U_\al$. 
If $\xi$ is purely excessive then, letting
$\al \longrightarrow 0$, we deduce that
$\xi'\circ(I+ \beta U)=\xi$. 
Let   $\eta\in\ex$ be such
that $\xi'\leq$ $\eta - \eta\circ \beta U_\beta$. The measure
$\eta_1=\eta- \inf_\al \eta\circ \al U_\al$  is purely excessive
and clearly $\eta - \eta\circ \beta U_\beta=$ $\eta_1 - \eta_1
\circ \beta U_\beta$. Therefore $\xi= \xi'\circ (I+\beta U)\leq$
$(\eta_1 - \eta_1\circ \beta U_\beta)\circ (I+\beta U)=\eta_1\leq
\eta.$

$b)$ Assume that  the measure $\xi= \xi'\circ (I+\beta U)$ is $\sigma$-finite and let $\al>0$.
Then $\xi\circ \al U_\al =$ $\xi'\circ (I+\beta U)\al U_\al=$
$\xi'\al \circ U_\al + \beta \xi'\circ (U- U_\al)$.
Therefore if $\al >\beta$ then
$\xi\circ \al U_\al=$ $\xi'\circ (\al - \beta)U_\al +\beta\xi'\circ U \leq$
$\xi'+ \beta\xi'\circ U=\xi$.
If $\al \leq \beta$ then  $\xi\circ \al U_\al\leq \beta\xi'\circ U\leq \xi'\circ (I+\beta U)=\xi$.
Consequently the measure $\xi$ is $\cu$-excessive.
From $\xi'\circ U_\al\leq \xi'\circ U$ we get
$\xi\circ \al U_\al\leq \al\xi'\circ U +\beta \xi'\circ (U-U_\al).$
The measure $\xi'\circ U$ being $\sigma$-finite we conclude that
$\inf_\al \xi\circ \al U_\al =0$.\\

The following proposition is close to the results of R. K. Getoor
from \cite{Ge 80} and \cite{Ge 90}.\\

\noindent
{\bf Proposition A1}. {\it
If $\xi \in \mbox{\sf Exc}_{\mathcal{U}}$ then the following assertions are equivalent.

$1)$ The measure $\xi$ is dissipative.

$2)$ If $f \in p{\cal B}$, $ f > 0$ on $E$ and $\xi (f) <\infty$
then $U f < \infty$  $\xi$-a.e.

$3)$  There exists $F \in \mathcal{B}$ such that $\xi (E\setminus
F)=0$, $U (1_{E\setminus F}) = 0$ on $F$ and ${\cal U}|_{F}$ is
proper.

$4)$  There exists a finely continuous function $f \in bp
\mathcal{B}$, $ f > 0$ on $E$ such that $Uf \leq 1$ $\xi$-a.e.

$5)$  There exists $f \in bp \mathcal{B}$ such that $U f > 0$ on
$E$ and $U f \leq 1$ $\xi$-a.e.

$6)$  There exists a sequence $(f_n)_n \subset p\mathcal{B}$ such
that $U f_n$ is  bounded $\xi$-a.e.  for all $n$ and $U f_n
\nearrow \infty$. }

\begin{proof}
The equivalence  $1) \Longleftrightarrow 2)$  follows  from (2.11)
in \cite{Ge 90}. The  implications $ 4) \Longrightarrow 5)
\Longrightarrow 6)$  are clear. We have $3) \Longrightarrow 1)$
since $\mbox{\sf Exc}_{\mathcal{U}} = \mbox{\sf
Diss}_{\mathcal{U}}$ if $\mathcal{U}$ is proper.

$2) \Longrightarrow 3).$  Let $ g \in bp \mathcal{B}$, $g > 0$ on
$E$ be such that $\xi(g) < \infty$. Then $U g  < \infty$
$\xi$-a.e. If we set $A_n = [Ug \leq n]$ then $(A_n)_n\subset
{\cal B}$ is an increasing  sequence, $\xi (E \setminus \cup_n
A_n) = 0$ and $U(g 1_{A_n} ) \leq n$ for all $n$. The function $f
= g (1_{A_\infty} +$ $\sum_{n \geq 1} \dfrac{1}{n 2^n} 1_{A_n} )$
is strictly positive and $U f \leq 1 $  on $[U g < \infty]$, where
$A_\infty = [Ug = \infty]$. Taking $E_o = [Uf \leq 1]$ and
applying Lemma 2.1 we obtain the required set $F$.

$3) \Longrightarrow 4).$  Let $g \in bp \mathcal{B}$, $ g >0$ on
$E$ be such that $Ug \leq 1$ $\xi$-a.e. The function $f = U_1g$ is
bounded, finely continuous, strictly positive and we have
$\xi$-a.e. $Uf= UU_1 g =$ $Ug - U_1 g \leq$ $Ug \leq  1.$

$6) \Longrightarrow 5).$ Let $(f_n)_n \subset p \mathcal{B}$ and
$(\alpha_n)_n \subset \mathbb{R}_+^*$ such that $U f_n \leq
\alpha_n$ $\xi$-a.e. for all $n$ and $Uf_n \nearrow \infty$.
Consider the function $f = \sum_n \dfrac{1}{\alpha_n 2^n} f_n$.
Clearly $f \in p \mathcal{B}$,  $Uf \leq 1$  $\xi$-a.e. and $Uf>0$
on $E$.

$5) \Longrightarrow 3).$ Let $g \in p \mathcal{B}$, $g \leq 1$, be
such that $Ug >0$ and $Ug \leq 1$ $\xi$-a.e., and let $F= [Ug <
\infty]$. 
From $Ug =U_\alpha g + \alpha U_\alpha Ug$ we get that
on $F$ we have $U_\alpha (1_{E \setminus F}) = 0$ and therefore
$U(1_{E \setminus F}) = 0$. The function $f=\alpha U_\alpha g
\cdot 1_F + 1_{E \setminus F}$ belongs to $p \mathcal{B}$, $f \leq
1$ and 
$Uf \leq \alpha U_\alpha Ug + U(1_{E\setminus F}) \leq Ug<\infty $
on $F$. It remains to show that $f > 0$. If we assume that $f(x)=
0$ then $x \in F$ and $U_\alpha g(x) = 0$. 
Consequently, we get $\alpha U_\alpha Ug(x) = Ug(x)$ 
and thus $\beta U_\beta Ug(x) =
Ug(x)$ for all $\beta >0$, $U_\beta g(x) = 0$. 
This leads to the contradictory equality $Ug(x) = 0$.
\end{proof}


\end{document}